\documentclass[number,citesort,seceqn,dvips]{arxbj}

\aid{0}
\volume{17}
\issue{3}
\pubyear{2011}
\firstpage{1054}
\lastpage{1062}
\doi{10.3150/10-BEJ304}

\makeatletter
\newremark{example}{Example}
\newremark{remark}{Remark}
\newtheorem{theorem}{Theorem}
\newtheorem{corollary}{Corollary}
\makeatother

\begin{document}
\begin{frontmatter}

\title{A note on a maximal Bernstein inequality}
\runtitle{A note on a maximal Bernstein inequality}

\begin{aug}
\author[1]{\fnms{P\'eter} \snm{Kevei}\thanksref{1}\corref{}\ead[label=e1,text=kevei@cimat.mx]{kevei@cimat.mx}} \and
\author[2]{\fnms{David M.} \snm{Mason}\thanksref{2}\ead[label=e2]{davidm@udel.edu}}
\runauthor{P. Kevei and D.M. Mason}
\address[1]{CIMAT, Callej\'on Jalisco S/N, Mineral de Valenciana,
Guanajuato 36240, Mexico \\ \printead{e1}}
\address[2]{Statistics Program, University of Delaware, 213 Townsend Hall
Newark, DE 19716, USA \\ \printead{e2}}
\end{aug}
\received{\smonth{2} \syear{2010}}
\revised{\smonth{6} \syear{2010}}

\vspace*{3.5pt}
\textit{Dedicated to the memory of S\'{a}ndor Cs\"{o}rg\H{o}}
\vspace*{-10pt}
\begin{abstract}
We show somewhat unexpectedly that whenever a general
Bernstein-type maximal inequality holds for partial sums of a sequence of random
variables, a maximal form of the inequality is also valid.
\end{abstract}


\end{frontmatter}

\section{Introduction and statement of main result}

Let $X_{1},X_{2},\dots,$ be a sequence of independent random variables
such that for all $i\geq1$, $EX_{i}=0$ and for some $\kappa>0$ and
$v>0$ for integers $m\geq2$, $E\vert X_{i}\vert ^{m}\leq
vm!\kappa^{m-2}/2$. The classic Bernstein inequality (cf. \cite{SW}, page 855) says that, in this situation, for all $n\geq1$ and $%
t\geq0,$%
\[
\mathbf{P}\Biggl\{ \Biggl\vert \sum_{i=1}^{n}X_{i}\Biggr\vert
>t\Biggr\} \leq2\exp\biggl\{ -\frac{t^{2}}{2vn+2\kappa t}\biggr\} .
\]
Moreover (cf. \cite{R}, Theorem B.2), its maximal form also holds; that
is, we have
\[
\mathbf{P}\Biggl\{ \max_{1\leq j\leq n}\Biggl\vert
\sum_{i=1}^{j}X_{i}\Biggr\vert >t\Biggr\} \leq2\exp\biggl\{ -\frac{t^{2}}{%
2vn+2\kappa t}\biggr\} .
\]
It turns out that, under a variety of assumptions, a sequence of not
necessarily independent random variables $X_{1},X_{2},\dots,$ will
satisfy a generalized Bernstein-type inequality of the following form:
For suitable constants $A>0$, $a>0$, $b\geq0$ and $0<\gamma<2$ for all
$m\geq0$, $n\geq1$ and $t\geq0$,
%
\begin{equation} \label{assumpgamma}
\mathbf{P}\{|S(m+1,m+n)|>t\}\leq A\exp\biggl\{ -\frac{at^{2}}{n+bt^{\gamma}}%
\biggr\} ,
\end{equation}
where, for any choice of $1\leq k\leq l<\infty$, we denote the partial sum $%
S(k,l)=\sum_{i=k}^{l}X_{i}.$ Here are some examples.

\begin{example}\label{ex1}
Let $X_{1},X_{2},\dots,$ be a stationary
sequence satisfying $EX_{1}=0$ and \mbox{$\operatorname{Var}X_{1}=1$}. For each integer
$n\geq1$ set $S_{n}=X_{1}+\cdots+X_{n}$ and $B_{n}^{2}=\operatorname{Var}( S_{n}) $.
Assume that for some $\sigma_{0}^{2}>0$ we have
$B_{n}^{2}\geq\sigma_{0}^{2}n $ for all $n\geq1$. Statulevi\v{c}ius and
Jakimavi\v{c}ius \cite{SJ} and
Saulis and Statulevi\v{c}ius~\cite{SS} prove that the partial sums satisfy (%
\ref{assumpgamma}) with constants depending on a Bernstein condition on
the moments of $X_{1}$ and the particular mixing condition that the
sequence may fulfill. In fact, all values of $1\leq\gamma<2$ are
attainable. Their Bernstein-type inequalities are derived via a result
of \cite{BR} relating cumulant behavior to tail
behavior, which says that for an arbitrary random variable $\xi$ with
expectation~$0$, whenever there
exist $\gamma\geq0$, $H>0$ and $\Delta>0$ such that its cumulants $%
\Gamma_{k}( \xi) $ satisfy $\vert \Gamma_{k}( \xi) \vert \leq( k!/2)
^{1+\gamma}H/\Delta^{k-2}$ for $k=2,3,\dots,$ then for all $x\geq0$
%
\begin{equation}\label{BR}
\mathbf{P}\{ \pm\xi>x\} \leq\exp\biggl\{ -\frac{x^{2}}{2( H+( x/\Delta^{1/(
1+2\gamma) }) ^{( 1+2\gamma) /(1+\gamma)}) }\biggr\} .
\end{equation}
\end{example}

In Example \ref{ex1}, $\xi=S_{n}/B_{n}$ and $\Delta=d\sqrt{n}$ for some $d>0$.

\begin{example}
Doukhan and Neumann \cite{DN} have shown,
using
the result in (\ref{BR}), that if a sequence of mean zero random variables $%
X_{1},X_{2},\dots,$ satisfies a general covariance condition, then the
partial sums satisfy (\ref{assumpgamma}). Refer to their Theorem 1 and
Remark 2, and also see \cite{KN}.
\end{example}

\begin{example}
Assume that $X_{1},X_{2},\dots,$ is a
strong mixing sequence with mixing coefficients $\alpha( n) $,
$n\geq1$, satisfying for some $c>0$, $\alpha( n) \leq\exp( -2cn)
$. Also assume that $EX_{i}=0$ and for some $M>0$ for all $i\geq1$, $%
\vert X_{i}\vert \leq M$. Theorem 2 of Merlev\'{e}de, Peligrad and Rio
\cite{MPR} implies that for some constant $C>0$ for all $t\geq0$ and
$n\geq1$,
%
\begin{equation} \label{exin2}
\mathbf{P}\{ |S_{n}|>t\} \leq\exp\biggl( -\frac{Ct^{2}}{%
nv^{2}+M^{2}+tM( \log n) ^{2}}\biggr) ,
\end{equation}
with $S_{n}=\sum_{i=1}^{n}X_{i}$ and where $v^{2}=\sup_{i>0}( \operatorname{Var}(
X_{i}) +2\sum_{j>i}\vert \operatorname{cov}( X_{i},X_{j}) \vert )
>0.$
\end{example}

To see how the last example satisfies (\ref{assumpgamma}), notice that
for any $0<\eta<1$ there exists a $D_{1}>0$ such that for all $t\geq0$
and $n\geq1$,
%
\begin{equation}\label{ex2}
nv^{2}+M^{2}+tM( \log n) ^{2}\leq n( v^{2}+M^{2}) +D_{1}t^{1+\eta}.
\end{equation}
Thus (\ref{assumpgamma}) holds with $\gamma=1+\eta$ for suitable $A>0$,
$a>0$ and $b\geq0$.

For any choice of $1\leq i\leq j<\infty$ define%
\[
M(i,j)=\max\{|S(i,i)|,\ldots,|S(i,j)|\}.
\]
We shall show, somewhat unexpectedly, that if a sequence of random
variables
$X_{1},X_{2},\ldots,$ satisfies a Bernstein-type inequality of the form (\ref%
{assumpgamma}), then, without any additional assumptions, a modified
version of it also holds for $M(1+m,n+m)$ for all $m\geq0$ and
$n\geq1$.

\begin{theorem}\label{T1}
 Assume that, for constants $A>0$, $a>0$, $b\geq0$ and
$\gamma\in(0,2)$, inequality~(\ref{assumpgamma}) holds for all $m\geq0,n\geq1$ and $%
t\geq0$. Then for every $0<c<a$ there exists a~$C>0$ depending only on $A,a$, $b$ and $\gamma$ such that for all $n\geq1$, $m\geq0$ and $t\geq0$,
%
\begin{equation} \label{Cineq}
\mathbf{P}\{M(m+1,m+n)>t\}\leq C\exp\biggl\{ -\frac{ct^{2}}{n+bt^{\gamma}}\biggr\} .
\end{equation}
\end{theorem}

\begin{remark}\label{RA}
Notice that though $c<a$,  $c$ can be chosen arbitrarily
close to $a$.
\end{remark}

\begin{remark}\label{R1}
Theorem \ref{T1} was motivated by Theorem 2.2 of M\'{o}ricz,
Serfling and Stout \cite{MSS}, who showed that whenever for a suitable
positive function $g( i,j) $ of $( i,j) \in \{
1,2,\dots \} \times \{ 1,2,\dots \} $, positive function $%
\phi ( t) $ defined on $( 0,\infty ) $ and constant $%
K>0$, for all $1\leq i\leq\break j<\infty $ and $t>0$,
\[
\mathbf{P}\{|S(i,j)|>t\}\leq K\exp \{ -\phi ( t) /g( i,j) \} ,
\]%
then there exist constants $0<c<1$ and $C>0$ such that for all $n\geq
1$ and $t>0$,
\[
\mathbf{P}\{M(1,n)>t\}\leq C\exp \{ -c\phi ( t) /g( 1,n) \} .
\]%
Earlier, M\'{o}ricz \cite{M79} proved that in the special case when
$\phi ( t) =t^{2}$ one can choose \mbox{$c<1$} arbitrarily close to~$1$ by
making $C>0$ large enough. This inequality is clearly not applicable to
obtain a maximal form of the generalized Bernstein inequality.
\end{remark}

\begin{remark}\label{R}
We do not know whether there exist examples for which (\ref%
{assumpgamma}) holds for some $0<\gamma<1$ and $b>0$. However, since
the proof of our theorem remains valid in this case, we shall keep it
in the statement.
\end{remark}

\begin{remark}\label{R2}
The version of Theorem \ref{T1} obtained by replacing everywhere $%
|S(m+1,m+n)|$ by $S(m+1,m+n)$ and $M(m+1,m+n)$ by $M^{+}( m+1,m+n)
=\max_{m+1\leq j\leq n+m}( S(m+1,j)\vee0) $ remains true with little
change in the proof.
\end{remark}

\begin{remark}\label{R3}
Theorem \ref{T1} also remains valid for sums of Banach space
valued random variables with absolute value $\vert \cdot\vert $
replaced by norm $\Vert \cdot\Vert $.
\end{remark}

\begin{remark}
\label{Statistics} In statistics, maximal exponential inequalities are
crucial tools to determine the exact rate of almost sure pointwise and
uniform consistency of kernel estimators of the density function and
the regression function. The literature in this area is huge. See, for
instance, \cite{DM92,DM94,EM00,EM05,GG,St}
and the references therein. These results only treat the case of i.i.d.
observations. Dependent versions of our maximal Bernstein inequalities
should be useful to determine exact rates of almost sure consistency of
kernel estimators based on data that possess a certain dependence
structure. In fact, some work in this direction has already been
accomplished in Section 4.2 of \cite{DN}. To carry
out such an application in the present paper is well beyond its scope.
\end{remark}

Theorem \ref{T1} leads to the following bounded law of the iterated
logarithm.

\begin{corollary}\label{C1}
Under the assumptions of Theorem \ref{T1}, with probability $1$,
%
\begin{equation}   \label{BLIL}
\limsup_{n\rightarrow\infty}\frac{|S(1,n)|}{\sqrt{n\log\log n}}\leq\frac {1}{\sqrt{a}}.
\end{equation}
\end{corollary}

\begin{remark}\label{R5}
In general, one cannot replace ``$\leq$'' by ``$=$'' in (\ref{BLIL}).
To see this, let $Y$, $Z_{1},Z_{2},\dots$ be a~sequence of
independent random variables such that $Y$ takes on the value~$0$ or
$1$ with probability~$1/2$ and $Z_{1},Z_{2},\dots$ are independent
standard normals. Now define $X_{i}=YZ_{i}$, $i=1,2,\ldots.$ It is
easily checked that assumption (\ref{assumpgamma}) is satisfied with $A=2,$ $a=1/2$, $b=0$ and $%
\gamma=1.$ When $Y=1$ the usual law of the iterated logarithm gives
with probability $1$,
%
\begin{equation}\label{ac}
\limsup_{n\rightarrow\infty}|S(1,n)|/\sqrt{n\log\log n}=\sqrt{2}=1/\sqrt {a},
\end{equation}
whereas, when $Y=0$ the above limsup is obviously $0.$ This agrees with
Corollary \ref{C1}, which says that with probability $1$ the limsup is
${\leq}\sqrt{2}$. However, we see that with probability $1/2$ it
equals~$\sqrt{2} $ and with probability $1/2$ it equals $0$.
\end{remark}

Theorem \ref{T1} is proved in Section \ref{sec2} and the proof of Corollary
\ref{C1} is given in Section \ref{sec3}.

\section{Proof of theorem}\label{sec2}

The case $b=0$ is a special case of Theorem 1 of~\cite{M79}. Therefore
we shall always assume that $b>0$. Choose any $0<c<a.$ We prove our
theorem by induction on $n$. Notice that by the assumption, for any
integer $n_{0}\geq1$ we may choose $C>An_{0}$ to make the statement
true for all $1\leq n\leq n_{0}$. This remark will be important,
because at some steps of the proof we
assume that $n$ is large enough. Also, since the constants $A$, $a$, $b$ and $%
\gamma$ in (\ref{assumpgamma}) are independent of $m$,  we can
 assume $m=0$ without loss of generality in our proof.

Assume the statement holds up to some $n\geq2$. (The constant $C$ will
be determined in the course of the proof.)

\textit{Case} 1: Fix a $t > 0$ for which
%
\begin{equation}  \label{alpha}
t^{\gamma} \leq\alpha n
\end{equation}
for some $0<\alpha<1$ to be specified later. (In any case, we assume that $%
\alpha n\geq1$.) Using an idea of~\cite{MSS}, we may write for arbitrary $%
1\leq k\leq n$, $0\leq q\leq1$ and $p+q=1$ the inequality
\begin{eqnarray*}
&&\mathbf{P}\{M(1,n+1)>t\}\\
&&\quad\leq\mathbf{P}\{M(1,k)>t\}+\mathbf{P}\{|S(1,k)|>pt\}+\mathbf{P}%
\{M(k+1,n+1)>qt\}.
\end{eqnarray*}

Let
\[
u=\frac{ n + t^{\gamma} b (q^{\gamma}-q^{2})}{1+q^{2}}.
\]
Notice that
%
\begin{equation} \label{eq}
\frac{t^{2}}{u+ b t^{\gamma}}=\frac{q^{2}t^{2}}{n-u + b q^{\gamma}t^{\gamma}}.
\end{equation}
Set
%
\begin{equation} \label{k}
k=\lceil u\rceil .
\end{equation}
Using the induction hypothesis and (\ref{assumpgamma}) we obtain
\begin{eqnarray}\label{mainineqgamma}
&&\mathbf{P}\{M(1,n+1)>t\}\nonumber\\ [-8pt]\\ [-8pt]
&&\quad\leq C\exp\biggl\{ -\frac{ct^{2}}{k+ b t^{\gamma}}\biggr\} + A \exp\biggl\{ -%
\frac{a p^{2}t^{2}}{k + b p^{\gamma}t^{\gamma}}\biggr\} +C\exp\biggl\{ -\frac{%
c q^{2} t^{2}}{n - k + b q^{\gamma}t^{\gamma}}\biggr\} .\nonumber
\end{eqnarray}
Notice that we chose $k$ to make the first and third terms in the right-hand side of (\ref%
{mainineqgamma}) almost equal, and since by (\ref{k})
\[
\frac{t^{2}}{k+ b t^{\gamma}}\leq\frac{q^{2}t^{2}}{n - k + b q^{\gamma
}t^{\gamma}},
\]
the first term is greater than or equal to the third.

First we handle the second term in (\ref{mainineqgamma}), showing that for $%
0\leq t\leq(\alpha n)^{1/\gamma}$,
\[
\exp\biggl\{ -\frac{ap^{2}t^{2}}{k+bp^{\gamma}t^{\gamma}}\biggr\} \leq \exp\biggl\{
-\frac{ct^{2}}{n+1+bt^{\gamma}}\biggr\}.
\]
For this we need to verify that for $0\leq t\leq(\alpha n)^{1/\gamma}$,
%
\begin{equation}\label{eq1}
\frac{ap^{2}}{k+bp^{\gamma}t^{\gamma}}>\frac{c}{n+1+bt^{\gamma}},
\end{equation}
which is equivalent to
\[
ap^{2}(n+1+bt^{\gamma})>c(k+bp^{\gamma}t^{\gamma}).
\]
Using that
\[
k=\lceil u\rceil\leq u+1=1+\frac{1}{1+q^{2}}[ n+b(q^{\gamma}-q^{2})t^{%
\gamma}] ,
\]
it is enough to show
\[
n\biggl( ap^{2}-\frac{c}{1+q^{2}}\biggr) +t^{\gamma}\biggl(
ap^{2}b-cbp^{\gamma}-\frac{cb}{1+q^{2}}(q^{\gamma}-q^{2})\biggr) +ap^{2}-c>0.
\]
Note that if the coefficient of $n$ is positive, then we can choose
$\alpha$ in (\ref{alpha}) small enough to make the above inequality
hold, even if the
coefficient of $t^{\gamma}$ is negative. So in order to guarantee (\ref{eq1}) (at least for large $n$) we only have to choose the parameter $p$ so
that $ap^{2}-c>0$ -- which implies that
%
\begin{equation}\label{assump-c-1}
ap^{2}-\frac{c}{1+q^{2}}>0
\end{equation}
holds -- and then select $\alpha$ small enough.

Next we treat the first and third terms in (\ref{mainineqgamma}).
By the remark above, it is enough to handle the first term. Let
us examine the
ratio of $C\exp\{-ct^{2}/(k+bt^{\gamma})\}$ and $C\exp\{-ct^{2}/(n+1+bt^{%
\gamma})\}$. Notice again that since $u+1\geq k$,
\begin{eqnarray*}
n+1-k & \geq& n-u=n-\frac{n+b(q^{\gamma}-q^{2})t^{\gamma}}{1+q^{2}} \\
& =& \frac{q^{2}n-b(q^{\gamma}-q^{2})t^{\gamma}}{1+q^{2}} \\
&\geq& n\frac{q^{2}-\alpha b(q^{\gamma}-q^{2})}{1+q^{2}} \\
& =&\!: c_{1}n.
\end{eqnarray*}
At this point we need that $0<c_{1}<1$. Thus we choose $\alpha$ small
enough so that
%
\begin{equation}
q^{2}-\alpha b(q^{\gamma}-q^{2})>0.   \label{assump-c-2}
\end{equation}

Also, using $t\leq(\alpha n)^{1/\gamma}$, we get the bound
\[
(n + 1 + b t^{\gamma}) ( k + b t^{\gamma}) \leq n^{2} ( 1 + \alpha
b)^{2} =: c_{2} n^{2},
\]
which holds if $n$ is large enough. Therefore, we obtain for the ratio
\[
\exp\biggl\{ -ct^{2}\biggl( \frac{1}{k+b t^{\gamma}}-\frac{1}{n+1+b t^{\gamma }%
}\biggr) \biggr\} \leq\exp\biggl\{ -\frac{c c_{1} t^{2}}{c_{2} n}\biggr\} \leq%
\mathrm{e}^{-1},
\]
whenever $c c_{1} t^{2}/(c_{2} n) \geq1$, that is, $t\geq\sqrt{c_{2} n/ (c c_{1})}$. Substituting back into (\ref{mainineqgamma}), for $t\geq\sqrt{%
c_{2} n / (c c_{1})}$ and $t\leq(\alpha n)^{1/\gamma}$ we obtain
\begin{eqnarray*}
&&\mathbf{P}\{M(1,n+1)>t\}\\
&&\quad\leq\biggl( \frac{2}{\mathrm{e}}C + A \biggr) \exp\{-ct^{2}/(n+1+ b
t^{\gamma})\}\leq C\exp\{-ct^{2}/(n+1+ b t^{\gamma})\},
\end{eqnarray*}
where the last inequality holds for $C> A \mathrm{e}/(\mathrm{e}-2)$.

Next assume that $t<\sqrt{c_{2}n/(cc_{1})}$. In this case, choosing $C$
large enough, we can make the bound $>1$, namely
\[
C\exp\biggl\{ -\frac{ct^{2}}{n+1+bt^{\gamma}}\biggr\} \geq C\exp\biggl\{ -%
\frac{cc_{2}n}{cc_{1}n}\biggr\} =C\mathrm{e}^{-c_{2}/c_{1}}\geq1,
\]
if $C>\mathrm{e}^{c_{2}/c_{1}}$.

\textit{Case} 2: Now we must handle the case $t>(\alpha
n)^{1/\gamma }$. Here we apply the inequality
\[
\mathbf{P}\{M(1,n+1)>t\}\leq\mathbf{P}\{M(1,n)>t\}+\mathbf{P}%
\{|S(1,n+1)|>t\}.
\]
Using assumption (\ref{assumpgamma}) and the induction hypothesis, we
have
\[
\mathbf{P}\{M(1,n+1)>t\}\leq C\exp\biggl\{ -\frac{ct^{2}}{n+bt^{\gamma}}%
\biggr\} +A\exp\biggl\{ -\frac{at^{2}}{n+1+bt^{\gamma}}\biggr\} .
\]
We will show that the right-hand side $\leq
C\exp\{-ct^{2}/(n+1+bt^{\gamma})\}$. For this it is enough to prove
\begin{equation} \label{less}
\exp\biggl\{ -ct^{2}\biggl( \frac{1}{n+bt^{\gamma}}-\frac{1}{n+1+bt^{\gamma}}\biggr) \biggr\}
+\frac{A}{C}\exp\biggl\{ -\frac{t^{2}(a-c)}{n+1+bt^{\gamma }}\biggr\} \leq1.
\end{equation}

First assume that $\gamma\leq1$. Using the bound following from
$t>(\alpha n)^{1/\gamma}$, we get
\[
\frac{t^{2}}{(n+bt^{\gamma})(n+bt^{\gamma}+1)}\geq\frac{t^{2}}{(\alpha
^{-1}+b)(2\alpha^{-1}+b)t^{2\gamma}}=:t^{2-2\gamma}c_{3}\geq c_{3}.
\]
We have that the right-hand side of (\ref{less}) for $a\geq c$ is less than
\[
\mathrm{e}^{-cc_{3}}+\frac{A}{C}\leq1
\]
for $C$ large enough.

For $1<\gamma<2$ we have to use a different argument. For $t$ large
enough (i.e., for $n$ large enough, since $t>(\alpha n)^{1/\gamma}$) we
have
\[
\exp\biggl\{ -\frac{ct^{2}}{(n+bt^{\gamma})(n+bt^{\gamma}+1)}\biggr\}
\leq\exp\{ -cc_{3}t^{2-2\gamma}\} \leq1-\frac{cc_{3}t^{2-2\gamma }}{2}.
\]
We also have for $C>A$,
\[
\frac{A}{C}\exp\biggl\{ -\frac{t^{2}(a-c)}{n+1+bt^{\gamma}}\biggr\} \leq\exp\biggl\{
-t^{2-\gamma}\frac{a-c}{2\alpha^{-1}+b}\biggr\} .
\]
It is clear that since $a>c$, for $t$ large enough, that is, for $n$
large enough,
\[
\frac{cc_{3}t^{2-2\gamma}}{2}>\exp\biggl\{ -t^{2-\gamma}\frac{a-c}{%
2\alpha^{-1}+b}\biggr\} .
\]
The proof is complete.

\section{Proof of corollary}\label{sec3}

Choose any $\lambda>1$ and set $m_{r}=\lceil \lambda^{r}\rceil $ for
$r=1,2,\ldots.$ Now, using inequality (\ref{Cineq}), we get
\begin{eqnarray*}
&&\mathbf{P}\bigl\{ M(1,m_{r}) >\sqrt{c^{-1}m_{r+1}\log\log m_{r}}\bigr\}\\
&&\quad\leq C \exp\biggl\{ -\frac{m_{r+1}\log\log m_{r}}{m_{r}+b(
c^{-1}m_{r+1}\log\log m_{r}) ^{\gamma/2}}\biggr\} .
\end{eqnarray*}
Since as $r\rightarrow\infty$
\[
\frac{m_{r+1}\log\log m_{r}}{m_{r}+b( c^{-1}m_{r+1}\log\log m_{r})
^{\gamma/2}}=\bigl( 1+\mathrm{o}( 1) \bigr) \lambda\log r,
\]
it is readily checked that for $r_{0}$ large enough so that $\log\log
m_{r_{0}}>0,$
\[
\sum_{r=r_{0}}^{\infty}\mathbf{P}\bigl\{ M(1,m_{r}) >\sqrt{%
c^{-1}m_{r+1}\log\log m_{r}}\bigr\} <\infty
\]
and thus, since $m_{r+1}/m_{r}=\lambda+\mathrm{o}( 1) $, we get by the
Borel--Cantelli lemma that with probability~$1$
%
\begin{equation}\label{aa}
\limsup_{r\rightarrow\infty}\frac{ M(1,m_{r})}{\sqrt{m_{r}\log\log m_{r}}}%
\leq\sqrt{\lambda c^{-1}}.
\end{equation}
Next we see that for all $r\geq r_{0}$
\[
\max_{m_{r}\leq n<m_{r+1}}\frac{|S(1,n)|}{\sqrt{n\log\log n}}\leq \frac{%
M(1,m_{r+1})}{\sqrt{m_{r}\log\log m_{r}}}.
\]
Thus by (\ref{aa}), with probability $1,$
\begin{eqnarray*}
&&\limsup_{r\rightarrow\infty}\max_{m_{r}\leq n<m_{r+1}}\frac{|S(1,n)|}{\sqrt{%
n\log\log n}}\\
&&\quad\leq\limsup_{r\rightarrow\infty}\frac{M(1,m_{r+1})}{\sqrt{m_{r}\log\log
m_{r}}}\\
&&\quad=\limsup_{r\rightarrow\infty}\frac{M(1,m_{r+1})}{\sqrt{m_{r+1}\log\log
m_{r+1}}}\frac{\sqrt{m_{r+1}\log\log m_{r+1}}}{\sqrt{m_{r}\log\log m_{r}}}%
\leq\lambda\sqrt{c^{-1}}.
\end{eqnarray*}
Hence, since $\lambda>1$ can be chosen arbitrarily close to $1$ and
$c<a$ arbitrarily close to $a$, we have proved (\ref{BLIL}).

\section*{Acknowledgements}

The authors are  grateful to the referee for a
number of penetrating comments and suggestions that greatly improved
the paper. Kevei's research was partially supported by the Analysis and
Stochastics Research Group of the Hungarian Academy of Sciences and
Mason's by an NSF grant.

\printhistory


\begin{thebibliography}{16}
\bibitem{BR}
Bentkus, R. and Rudzkis, R. (1980). On exponential estimates of
the distribution of random variables. \textit{Lithuanian Math. J.}
\textbf{20} 15--30 (in Russian).
\MR{0575427}

\bibitem{DM92}
Deheuvels, P. and Mason, D.M. (1992). Functional laws of the
iterated logarithm for the increments of empirical and quantile processes.
\textit{Ann. Probab}. \textbf{20} 1248--1287.
\MR{1175262}

\bibitem{DM94}
Deheuvels, P. and Mason, D.M. (1994). Functional laws of the
iterated logarithm for local empirical processes indexed by sets.
\textit{Ann. Probab}. \textbf{22} 1619--1661.
\MR{1303659}

\bibitem{DN}
Doukhan, P. and Neumann, M.H. (2007). Probability and moment
inequalities for sums of weakly dependent random variables, with
applications. \textit{Stochastic Process. Appl}. \textbf{117} 878--903.
\MR{2330724}

\bibitem{EM00}
Einmahl, U. and Mason, D.M. (2000). An empirical process
approach to the uniform consistency of kernel-type function estimators.
\textit{J. Theoret. Probab.} \textbf{13} 1--37.
\MR{1744994}

\bibitem{EM05}
Einmahl, U. and Mason, D.M. (2005). Uniform in bandwidth
consistency of kernel-type function estimators. \textit{Ann. Statist.}
\textbf{33} 1380--1403.
\MR{2195639}

\bibitem{GG}
Gin\'{e}, E. and Guillou, A. (2002).
Rates of strong uniform
consistency for multivariate kernel density estimators. En l'honneur de J.
Bretagnolle, D. Dacunha-Castelle, I. Ibragimov.
\textit{Ann. Inst. H. Poincar\'{e} Probab. Statist}. \textbf{38} 907--921.
\MR{1955344}

\bibitem{KN}
Kallabis, R. and Neumann, M.H. (2006). An exponential
inequality under weak dependence. \textit{Bernoulli} \textbf{12} 333--350.
\MR{2218558}

\bibitem{MPR}
Merlev\`{e}de, F., Peligrad, M. and Rio, E. (2009). Bernstein
inequality and moderate deviations under strong mixing conditions. In
\textit{High Dimensional Probability V: The Luminy Volume} (C. Houdr\'{e},
V.~Koltchinskii, D.M. Mason and M. Peligrad, eds.)  273--292. Beachwood, OH:
IMS.

\bibitem{M79}
M\'{o}ricz, F.A. (1979).
Exponential estimates for the
maximum of partial sums. I. Sequences of rv's. Special issue dedicated to
George Alexits on the occasion of his 80th birthday. \textit{Acta Math.
Acad. Sci. Hungar.} \textbf{33} 159--167.
\MR{0515130}

\bibitem{MSS}
M\'{o}ricz, F.A., Serfling, R.J. and Stout, W.F. (1982).
Moment and probability bounds with quasisuperadditive structure for the
maximum partial sum. \textit{Ann. Probab. }\textbf{10} 1032--1040.
\MR{0672303}

\bibitem{R}
Rio, E. (2000). \textit{Th\'{e}orie asymptotique des processus al\'{e}atoires faiblement d\'{e}pendants. (French) Math\'{e}matiques \&
Applications (Berlin)} \textbf{31}. Berlin: Springer.
\MR{2117923}

\bibitem{SW}
Shorack, G.R. and Wellner, J.A. (1986). \textit{Empirical
Processes with Applications to Statistics}. New York: Wiley.
\MR{0838963}

\bibitem{SS}
Saulis, L. and Statulevi\v{c}ius, V.A. (1991). \textit{Limit
Theorems for Large Deviations}. Dordrecht: Kluwer.
\MR{1171883}

\bibitem{SJ}
Statulevi\v{c}ius, V.A. and Jakimavi\v{c}ius, D.A. (1988).
Estimates for semiinvariants and centered moments of stochastic processes
with mixing. I. \textit{Litovsk. Mat. Sb.} \textbf{28} 112--129;
translation in \textit{Lithuanian Math.~J}. \textbf{28} 67--80.
\MR{0949647}

\bibitem{St}
Stute, W. (1982). A law of the logarithm for kernel density
estimators. \textit{Ann. Probab}. \textbf{10} 414--422.
\MR{0647513}
\end{thebibliography}
\end{document}